\documentclass[twocolumn]{autart}  
\usepackage{graphicx}   
\usepackage{subfigure}       
\usepackage{amsmath}
\usepackage{amssymb}
\usepackage{bm}
\usepackage{bbm}
\usepackage{mathrsfs}
\usepackage[sort&compress]{natbib}
\usepackage{fixltx2e}
\usepackage{xcolor}

\newcommand{\new}[1] {\textcolor{blue}{{#1}}}
\usepackage{hyperref}
\makeatletter
\newcommand{\onetagright}{\tagsleft@false}

\makeatother

\begin{document}

\begin{frontmatter}
  \title{Projected Primal-Dual Gradient Flow of Augmented Lagrangian with Application to Distributed Maximization of the Algebraic Connectivity of a Network\thanksref{footnoteinfo}} 

  \thanks[footnoteinfo]{This work is supported by China Scholarship Council.}

  \author[KTH_M]{Han Zhang}\ead{hanzhang@kth.se},
  \author[KTH_A]{Jieqiang Wei}\ead{jieqiang@kth.se},
  \author[Toronto]{Peng Yi}\ead{peng.yi@utoronto.ca},
  \author[KTH_M]{Xiaoming Hu}\ead{hu@kth.se},              

  \address[KTH_M]{Department of Mathematics, KTH Royal Institute of Technology, SE-100 44, Stockholm, Sweden}     
  \address[KTH_A]{Department of Automatic Control, KTH Royal Institute of Technology, SE-100 44, Stockholm, Sweden}      
  \address[Toronto]{Department of Electrical and Systems Engineering, Washington University in St. Louis, USA}


  \begin{abstract}
	In this paper, a projected primal-dual gradient flow of augmented Lagrangian is presented to solve convex optimization problems that are not necessarily strictly convex.
	The optimization variables are restricted by a convex set with computable projection operation on its tangent cone as well as equality constraints.
	As a supplement of the analysis in \citep{niederlander2016distributed}, we show that the projected dynamical system converges to one of the saddle points and hence finding an optimal solution.
	Moreover, the problem of distributedly maximizing the algebraic connectivity of an undirected network by optimizing the port gains of each nodes (base stations) is considered. The original semi-definite programming (SDP) problem is relaxed into a nonlinear programming (NP) problem that will be solved by the aforementioned projected dynamical system. Numerical examples show the convergence of the aforementioned algorithm to one of the optimal solutions. The effect of the relaxation is illustrated empirically with numerical examples. A methodology is presented so that the number of iterations needed to reach the equilibrium is suppressed. Complexity per iteration of the algorithm is illustrated with numerical examples.
  \end{abstract}

\end{frontmatter}

\section{Introduction}

When solving a convex minimization problem with strong duality, it is well-known that the optimal solution is the saddle point of the Lagrangian. Hence 
it is natural to consider the gradient flow of Lagrangians (also known as saddle point dynamics) where the primal variable follows the negative gradient flow while the dual variable follows the gradient flow.
Gradient flow of Lagrangians is first studied by \citep{arrow1959studies}, \citep{kose1956solutions} and has been revisited by \citep{feijer2010stability}.
\citep{feijer2010stability} studies the case of strictly convex problems and provides methodologies to transform non-strictly convex problems to strictly convex problems to fit the framework. The convergence is shown by employing the invariance principle for hybrid automata. \citep{cherukuri2016asymptotic} studies the same strictly convex problem from the perspective of projected dynamical systems and is able to show the convergence by a LaSalle-like invariant principle for Carath\'eodory solutions. Instead of considering discontinuous dynamics, \citep{durr2011smooth} proposes a smooth vector field for seeking the saddle points of strictly convex problems. \citep{wang2011control} considers a strictly convex problem with equality constraints and with inequality constraints respectively. Saddle point dynamics is also used therein, however, it is worth noticing that their problem is still strictly convex. When they consider the problem with inequality constraints, logarithmic barrier function is used. Though considering nonsmooth problems, \citep{zeng2017distributed} uses the projected saddle point dynamics of augmented Lagrangian whose equality constraint is the variable consensus constraint, and can be viewed as a special case of our problem.
Instead of using the continuous-time saddle point dynamics, an iterative distributed augmented Lagrangian method is developed in \citep{chatzipanagiotis2015augmented}. In a recent work \citep{niederlander2016distributed} and its conference version \citep{niederlander2016exponentially}, the authors consider the nonsmooth case of projected saddle point dynamics and the dynamics are the same as the ones in the current paper when the objective function is smooth. 

\vspace{-0.2cm}

In this paper, we will focus on maximizing network algebraic connectivity distributedly. In \citep{simonetto2013constrained}, the authors maximize the algebraic connectivity of a mobile robot network distributedly. The authors use first-order Taylor expansion to approximate the original non-convex problem and get a convex problem. A more general linear dynamics are considered and a two-step algorithm is proposed to solve the problem distributedly. It is shown in \citep{simonetto2013constrained} that the algebraic connectivity is monotonically increasing with the algorithm, while the convergence to one optimal solution is not explicitly given. \citep{schuresko2008distributed},\citep{yang2010decentralized} and \citep{zavlanos2008distributed} focus on assuring the connectivity distributedly, while the algebraic connectivity maximization is not considered.

\vspace{-0.2cm}

The main contribution of this paper is as follows. As a supplement to \citep{niederlander2016distributed} and its conference version \cite{niederlander2016exponentially}, we propose a novel analysis line regarding the convergence of the dynamical system to reach comparable results.
Moreover, the problem of distributedly maximizing the algebraic connectivity of an undirected network by adjusting the ``port gains" of each nodes (base stations) is considered. It is worth noticing that the problem motivates from a physical system and the goal is to enable each base station to compute its own optimal port gains only using its neighbours’
information, the total number of nodes $N$ and the information belonging to itself; one can not ``design'' the communication network according to the structure of the problem or the algorithm. (For example, \citep{pakazad2015distributed}).
We solve the original problem, which is an SDP, by relaxing it into an NP problem. The NP problem is not strictly convex, hence we adapt the projected saddle point dynamics method proposed in this work to solve the aforementioned NP problem. Numerical examples show that the aforementioned algorithm converges to one of the optimal solutions.

\section{Preliminaries and Notations}\label{sec:preliminaries}
We denote $\mathbbm{1}=\bm{1}\bm{1}^T$ as an $N$ dimensional all-one matrix, where $\bm{1}$ is an $N$ dimensional all-one vector. The element located on the $i$th row and and $j$th column of a matrix $A$ is denoted as $[A]_{ij}$. If matrix $A_1-A_2$ is positive semi-definite, then it will be denoted as $A_1\succeq A_2$. We use $\|\cdot\|$ to denote 2-norm of vectors. $|S|$ denotes the cardinality of set $S$. And any notation with the superscript $*$ is denoted as the optimal solution to the corresponding optimization problem. $tr(\cdot)$ denotes the trace of a matrix. $\langle\cdot,\cdot\rangle_2$ is denoted as the inner-product in Euclidean space and $\langle A_1,A_2\rangle_M=tr(A_1A_2)$ denotes the inner-product in $\mathcal{S}^n$, which is the Hilbert space of $n\times n$ symmetric matrix. 

\vspace{-0.2cm}

Assume $K\subset \mathbb{R}^n$ is a closed and convex set, the projection of a point $x$ to the set $K$ is defined as $P_K(x)=\arg\min_{y\in K}\|x-y\|$. For $x\in K$, $v\in \mathbb{R}^n$, the projection of the vector $v$ at $x$ with respect to $K$ is defined as: ( see \citep{nagurney2012projected},\citep{brogliato2006equivalence}) $\Pi_K(x,v)=\lim_{\delta\rightarrow 0}\frac{P_K(x+\delta v)-x}{\delta}=P_{T_K(x)}(v)$, where $T_K(x)$ denotes the tangent cone of $K$ at $x$. The interior, the boundary and the closure of $K$ is denoted as $int(K)$, $\partial K$ and $cl(K)$, respectively.
The set of inward normals of $K$ at $x$ is defined as $n(x)=\big\{\gamma\:|\:\|\gamma\|=1,\langle\gamma,x-y\rangle_2\leq 0,\forall y\in K \big\}$, and $\Pi_K(x,v)$ fulfills the following lemma:
\begin{lem}[\cite{nagurney2012projected}]\label{lem:Pi expression}
  If $x\in int(K)$, then $\Pi_K(x,v)=v$; if $x\in\partial K$, then $\Pi_K(x,v)=v+\beta(x)n^*(x)$, where $n^*(x)=\arg\max_{n\in n(x)}\langle v,-n\rangle$ and $\beta(x)=\max\{0,\langle v,-n^*(x)\rangle\}$.
\end{lem}

\vspace{-0.2cm}

Let $F$ be a vector field such that $F:K\mapsto \mathbb{R}^n$, the projected dynamical system is given by $\dot{x}=\Pi_K(x,F(x))$.
Note that the right hand side of above dynamics can be discontinuous on the $\partial K$. Hence given an initial value $x_0\in K$, the system does not necessarily have a classical solution. However, if $F(x)$ is Lipschitz continuous, then it has a unique Carath\'eodory solution that continuously depends on the initial value \citep{nagurney2012projected}.
\section{Problem Formulation and Projected Saddle Point Dynamics}\label{sec:problem_formulation}
In this section, we consider the following optimization problem defined on $\mathbb{R}^n$:
\begin{equation}
  \begin{aligned}
	& \underset{x\in K}{\text{minimize}}
	& & f(x)\\
	& \text{subject to}
	& & Ax-b=0,\\
  \end{aligned}
  \label{eq:classical_op}
\end{equation}
where $f:\mathbb{R}^n\mapsto\mathbb{R}$ and $A\in\mathbb{R}^{m\times n}$. $K$ is a convex set such that calculating the projection on its tangent cone is computationally cheap. $f(x)$ is a convex function but not necessarily strictly convex. It is also assumed that the gradient of $f(x)$ is locally Lipschitz continuous and the Slater's condition holds for \eqref{eq:classical_op}. Hence strong duality holds for \eqref{eq:classical_op}.

\vspace{-0.2cm}

The Lagrangian $\mathscr{L}:K\times \mathbb{R}^m\mapsto\mathbb{R}$ for the problem \eqref{eq:classical_op} is given by
\begin{equation}
  \mathscr{L}(x,v)=f(x)+v^T(Ax-b),
  \label{eq:lagrangian}
\end{equation} 
where $v\in\mathbb{R}^m$ is the Lagrangian multiplier of the constraint $Ax-b=0$. Since strong duality holds for \eqref{eq:classical_op}, then $(x^*,v^*)$ is a saddle point of $\mathscr{L}(x,v)$ if and only if $x^*$ is an optimal solution to \eqref{eq:classical_op} and $v^*$ is optimal solution to its dual problem. 
The augmented Lagrangian $\mathscr{L}_\mathcal{A}:K\times\mathbb{R}^m\mapsto\mathbb{R}$ for \eqref{eq:classical_op} is given by $\mathscr{L}_\mathcal{A}(x,v)=f(x)+v^T(Ax-b)+\frac{\rho}{2}(Ax-b)^T(Ax-b)$, where $\rho>0$ is the damping parameter that will help to suppress the oscillation of $x$ during optimization algorithms. Without loss of generality, we choose $\rho=1$.

\vspace{-0.2cm}

We propose to find the saddle point of \eqref{eq:lagrangian} via the saddle point dynamics projected on the set $K$, i.e.,
\begin{subequations}\label{eq:projected dyamics}
  \begin{align}
	\dot{x}&=\Pi_K(x,-\nabla f(x)-A^Tv-A^T(Ax-b))\nonumber\\
	&=\Pi_K(x,-\frac{\partial\mathscr{L}_\mathcal{A}(x,v)}{\partial x}),\\
	\dot{v}&=Ax-b=\frac{\partial\mathscr{L}_\mathcal{A}(x,v)}{\partial v}.
  \end{align}
\end{subequations}
Note that it is assumed that $\nabla f(x)$ is locally Lipschitz continuous, therefore there is a unique Carath\'eodory solution for the dynamics \eqref{eq:projected dyamics}.
\section{Convergence Analysis}\label{sec:convergence analysis}
In this section, we analyse the convergence for \eqref{eq:projected dyamics} and start with the analysis of the equilibrium point of \eqref{eq:projected dyamics}. \cite{niederlander2016distributed} consider the nonsmooth case of projected saddle point dynamics and the dynamics are the same as the ones in the current paper when the objective function is smooth. As a supplement, we propose a novel analysis line regarding the stability of the dynamical system to reach comparable results.
\begin{prop}\label{prop:opt sol equilibrium}
  $(x^*,v^*)$ is a saddle point to \eqref{eq:classical_op} if and only if it is an equilibrium of \eqref{eq:projected dyamics}.
\end{prop}
\begin{pf}
  Since strong duality holds for \eqref{eq:classical_op}, the optimality conditions become necessary and sufficient conditions. The optimality condition for \eqref{eq:classical_op} is given by $-\nabla f(x^*)-A^Tv^*\in N_K(x^*),\quad Ax^*-b=0$, \citep{eskelinen2007andrzej}, which implies $-\nabla f(x^*)-A^Tv^*+A^T(Ax^*-b)\in N_K(x^*)$, where $N_K(x^*)$ denotes the normal cone of $K$ at $x^*$.
  This implies $\Pi_K(x^*,-\nabla f(x^*)-A^Tv^*-A^T(Ax^*-b))=0$, therefore, $(x^*,v^*)$ is an equilibrium point of \eqref{eq:projected dyamics}. On the other hand, if $(x^*,v^*)$ is an equilibrium point of \eqref{eq:projected dyamics}, it must have $-\nabla f(x^*)-A^Tv^*+A^T(Ax^*-b)\in N_K(x^*)$ and $Ax^*-b=0$, which implies the optimality condition.
  \qed
\end{pf}
  \begin{prop}\label{thm:convergence}
	Given an initial value $(x(0),v(0))$, where $x(0)\in K$, the trajectory of the projected dynamical system \eqref{eq:projected dyamics} asymptotically converges to one of the saddle points of \eqref{eq:classical_op}.
  \end{prop}
  \begin{pf}
  We use LaSalle invariance principle for Carath\'eodory solutions \citep{cherukuri2015convergence} to prove the proposition.
	Suppose $(x^*,v^*)$ is a saddle point of the Lagrangian \eqref{eq:lagrangian}, namely, $x^*$ is the optimal solution of \eqref{eq:classical_op} and $v^*$ is the optimal solution of its dual problem. Construct the following Lyapunov function
	\begin{equation}
		d(x,v)=\frac{1}{2}(\|x-x^*\|^2+\|v-v^*\|^2).
	  \label{eq:lyap}
	\end{equation}
	Note that $d(x,v)$ is continuously differentiable and denote the right hand side of the dynamics \eqref{eq:projected dyamics} as vector field $F$. The Lie derivative along the vector field $f(x)$ of a function $V(x)$ is defined as $\mathcal{L}_fV(x)=\frac{\partial V(x)}{\partial x}f(x)$.
	By the definition of saddle points, $\mathscr{L}(x^*,v)\leq\mathscr{L}(x^*,v^*)\leq\mathscr{L}(x,v^*)$. The Lie derivative of $d(x,v)$ along the vector field $F$ is given by $\mathcal{L}_Fd(x,v)=(x-x^*)^T\Pi_K(x,-\nabla f(x)-A^Tv-A^T(Ax-b))+(v-v^*)^T(Ax-b)$.
	By Lemma \ref{lem:Pi expression}, it holds that $\Pi_K(x,-\nabla f(x)-A^Tv-A^T(Ax-b))=-\nabla f(x)-A^Tv-A^T(Ax-b)+\beta(x)n^*(x)$, where $\beta(x)\geq 0$ and $\langle n^*(x),x-y\rangle\leq 0,\:\forall y\in K$. Since $Ax^*-b=0$, hence it follows that $\mathcal{L}_Fd(x,v)=(x-x^*)^T(-\nabla f(x)-A^Tv-A^T(Ax-b)+\beta(x)n^*(x))+(v-v^*)^T(Ax-b)\leq (x-x^*)^T(-\nabla f(x)-A^Tv-A^T(Ax-b))+(v-v^*)^T(Ax-b)=-(x-x^*)^T\frac{\partial\mathscr{L}(x,v)}{\partial x}-(Ax-b)^T(Ax-b)+(v-v^*)^T \frac{\partial\mathscr{L}(x,v)}{\partial v}.$
	
	\vspace{-0.2cm}
	Since $\mathscr{L}(x,v)$ is convex with respect to $x$ and concave with respect to $v$, it follows from the first order property of convex and concave function \citep{boyd2004convex} that $\mathcal{L}_Fd(x,v)\leq \mathscr{L}(x^*,v)-\mathscr{L}(x,v^*)-(Ax-b)^T(Ax-b)$.
	Since $\mathscr{L}(x^*,v)\leq\mathscr{L}(x,v^*)$, we have $\mathcal{L}_Fd(x,v)\leq -(Ax-b)^T(Ax-b)\leq 0$. Note that $d(x,v)$ is convex and differentiable. By Proposition \ref{prop:opt sol equilibrium} and using the result in \citep{BACCIOTTI2004841}, we can conclude that any saddle point is Lyapunov stable. Recall that $f(x)$ has a locally Lipschitz continuous gradient, hence the uniqueness of the Carath\'eodory solution to \eqref{eq:projected dyamics} can be guaranteed.
	Since $d(x,v)$ is differentiable and by the definition of Carath\'eodory solution, $x(t),v(t)$ are absolute continuous, $d(x(t),v(t))$ is differentiable almost everywhere with respect to $t$ and $\frac{d}{dt}d(x(t),v(t))=\mathcal{L}_Fd(x,v)$ holds almost everywhere on $t\geq 0$. Therefore $d(x(t),v(t))$ is continuous and non-increasing with respect to time.
	Note that $d(x,v)$ is radially unbounded, hence the set $\hat{\mathcal{S}}=\{(x,v)\in\mathbb{R}^n\times\mathbb{R}^m|d(x,v)\leq d(x(0),v(0))\}$ is a compact invariant set for the system \eqref{eq:projected dyamics}. 
	The invariance of the $\omega$-limit set \citep{cherukuri2015convergence} can be proved by using the same methodology as Lemma 4.1 in \citep{khalil1996nonlinear} which is based on the continuity and the uniqueness of the solution. 
	Now by LaSalle invariance principle for Carath\'eodory solutions \citep{cherukuri2015convergence}, the trajectory of \eqref{eq:projected dyamics} converges to the largest invariant set in $cl(\Omega)$, where $\Omega=\{(x,v)\in\hat{\mathcal{S}}|\mathcal{L}_Fd(x,v)=0\}$.

\vspace{-0.2cm}
For $(x^\prime,v^\prime)\in \Omega$, we have $Ax^\prime-b=0$ and $\mathscr{L}(x^*,v^\prime)=\mathscr{L}(x^\prime,v^*)\Leftrightarrow f(x^*)=f(x^\prime)+v^{*T}(Ax^\prime-b)\Leftrightarrow f(x^*)=f(x^\prime)$.
Since $x^*$ is an optimal solution to \eqref{eq:classical_op}, $x^\prime$ is also an optimal solution to \eqref{eq:classical_op}.
	Denote $\mathcal{D}=\{(x,v)\in\hat{\mathcal{S}}|x \mbox{ is an optimizer}\}$ and note that $\Omega\subset \mathcal{D}$. Hence $cl(\Omega)\subset cl(\mathcal{D})$.
	On the other hand, since the set of optimal solutions for a convex optimization problem is closed, $\mathcal{D}$ is also closed and hence $cl(\Omega)\subset cl(\mathcal{D})=\mathcal{D}$. 

\vspace{-0.2cm}

	Denote $\mathcal{M}$ as the largest invariant set in $cl(\Omega)$. 
	Assume the initial value $(\hat{x}(0),\bar{v}(0))\in\mathcal{M}\subset cl(\Omega)$, then $\hat{x}(0)$ is also an optimal solution to \eqref{eq:classical_op}. Hence there exists some $\hat{v}$ such that $(\hat{x}(0),\hat{v})$ is a saddle point of \eqref{eq:lagrangian}. Then the Lyapunov function can be constructed similarly as $\hat{d}(x,v)=\frac{1}{2}(\|x-\hat{x}(0)\|^2+\|v-\hat{v}\|^2)$.
	We have shown that for any arbitrarily chosen saddle point $(x^*,v^*)$, $d(x(t),v(t))$ is non-increasing with respect to time; such statement also holds for $\hat{d}(x(t),v(t))$. 
	Furthermore, since $(\hat{x}(0),\bar{v}(0))\in \mathcal{M}$ and $\mathcal{M}$ is invariant, $(\hat{x}(t),\bar{v}(t))\in \mathcal{M},\forall t\geq 0$. This implies that $\hat{x}(t)$ is an optimizer for all $t\geq 0$ and hence we have $\dot{\bar{v}}(t)=A\hat{x}(t)-b=0,\forall t\geq 0$. And this implies $\bar{v}(t)=\bar{v}(0),\forall t\geq 0$ and hence $\hat{d}(\hat{x}(t),\bar{v}(t))\leq \hat{d}(\hat{x}(0),\bar{v}(0))=\frac{1}{2}\|\bar{v}(0)-\hat{v}\|^2$. This implies $\|\hat{x}(t)-\hat{x}(0)\|^2=0$ and hence $x(t)=\hat{x}(0),\forall t\geq 0$. Therefore, any trajectory that starts in $\mathcal{M}$ remains constant for all times, i.e., any point in $\mathcal{M}$ is an equilibrium point. By Proposition \ref{prop:opt sol equilibrium}, we can conclude that they are saddle points of \eqref{eq:classical_op}.

\vspace{-0.2cm}

	We have shown that given an initial value $(x(0),v(0))$, where $x(0)\in K$, the trajectory of \eqref{eq:projected dyamics} asymptotically converges to a set $\mathcal{M}$ whose elements are saddle points of \eqref{eq:classical_op}. Now we show the trajectory asymptotically converges to a point in $\mathcal{M}$ by contradiction. To abbreviate the notation, we denote $\eta=(x^T,v^T)^T$.
	Suppose $\eta(t)$ does not converge to a point in $\mathcal{M}$, namely, the trajectory's $\omega$-limit set $\Gamma(\eta(t))$ is not a singleton. This means we can choose $\bar{\eta}_1,\bar{\eta}_2\in\Gamma(\eta(t))\subset\mathcal{M}$, such that $\|\bar{\eta}_1-\bar{\eta}_2\|=\zeta,\zeta>0$. Since $\bar{\eta}_1,\bar{\eta}_2\in\mathcal{M}$, $\bar{\eta}_1,\bar{\eta}_2$ are saddle points and we have shown that all saddle points are Lyapunov stable.
	This means that there exists $\delta(\zeta/2)$, such that if $\|\eta(T)-\bar{\eta}_1\|<\delta(\zeta/2)$, then $\|\eta(t)-\bar{\eta}_1\|<\zeta/2,\forall t\geq T$. Since $\bar{\eta}_1$ is an $\omega$-limit point in $\Gamma(\eta(t))$, there exists such $T$ so that $\|\eta(T)-\bar{\eta}_1\|<\delta(\zeta/2)$, and hence the trajectory can never leave $\zeta/2$-neighbourhood of $\bar{\eta}_1$ after time instant $T$. But $\bar{\eta}_2$ is also an $\omega$-limit point, there must exists a sequence of point on the trajectory that tends to it. Hence we have a contradiction and $(x(t),v(t))$ converges to a saddle point in $\mathcal{M}$.
	\qed
  \end{pf}   
\section{Distributed Algebraic Connectivity Maximization}\label{sec:application}
  In this section, we apply the aforementioned algorithm to maximize the algebraic connectivity of a network in a distributed manner. The problem is first formulated as a Semi-definite Programming (SDP) problem. With an equivalent formulation of the original SDP problem, the problem is relaxed into a Nonlinear Programming (NP) problem in order to apply the aforementioned projected saddle point dynamics.
  \subsection{Motivation and Modeling}
\begin{figure}[!htpb]
	\centering
	\includegraphics[width=0.34\textwidth]{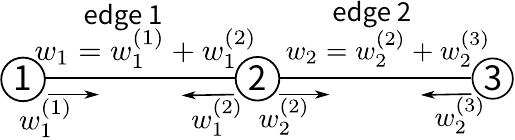}
	\caption{Edge weight is the sum of weights contributed by the nodes connected by the edge.}
	\label{fig:illustration}
\end{figure}

  The problem motivates from a physical communication network. Consider an undirected communication network $\mathcal{G(V,E,W)}$ whose nodes $i\in\mathcal{V}=\{1,2,\cdots,N\}$ are homogeneous base stations and can control their communication port gains $w_k^{(i)}\in\mathcal{W}$. (To abbreviate the notation, the edges are labelled with numbers.) The set of neighbours of node $i$ is denoted as $\mathcal{N}(i)$. The set of edges (communication channel) adjacent to node $i$ is denoted as $\mathcal{E}(i)$ and $\mathcal{E}=\bigcup_{i\in\mathcal{V}}\mathcal{E}(i)$. 
  As illustrated in Fig. \ref{fig:illustration}, the communication gain (strength) on each link $k\in\mathcal{E}$ is the sum of the port gains $w_k^{(i)}$ and $w_k^{(j)}$, $(i,j)=k\in\mathcal{E}$ contributed by the two end nodes connected by the edge. It is assumed that each agent can only get access to the information of its neighbours as well as the information of itself. Our goal is to develop a method so that each base station can adjust its own port gains only according to its neighbours' information, the number of nodes $N$ and the information belonging to itself, so that the algebraic connectivity of the total communication network is maximized.
  The graph $\mathcal{G(V,E,W)}$ we consider is undirected, and hence the weighted Laplacian matrix $L_w$ is symmetric and can be expressed as
  \begin{equation}\nonumber
	\begin{aligned}
	  \left[L_w\right]_{ij}&=
	  \begin{cases}
		\sum_lw_{il} &\mbox{if}\: i= j\:\mbox{and}\:(i,l)\in \mathcal{E}\\
		-w_{ij}&\mbox{if}\: i\neq j\:\mbox{and}\:(i,j)\in \mathcal{E}\\
		0 &\mbox{otherwise}
	  \end{cases}
	\end{aligned}
	\label{eq:Laplacian_matrix}
  \end{equation}
  and $L_w=\sum_{k\in \mathcal{E}}w_{k}E_{k}$, where $k$ is the label of the edges, and $0\leq w_k\in \mathcal{W},\:\forall k\in \mathcal{E}$ are the edge-weights. If node $i$ and $j$ are connected via edge $k$, then $[E_k]_{ii}=[E_k]_{jj}=1,\quad [E_k]_{ij}=[E_k]_{ji}=-1$, and the other elements of $E_k$ are zero. If the graph is connected, then the eigenvalues of $L_w$ satisfy: $0=\lambda_1<\lambda_2\leq\cdots\leq \lambda_N$ and $\lambda_2$ is the \textit{algebraic connectivity} of $\mathcal{G(V,E,W)}$. Let $L=\sum_{k\in\mathcal{E}}E_k$ be the unweighted Laplacian matrix of the graph. We suppose $L$ has only one zero eigenvalue, namely, the unweighted graph $\mathcal{G(V,E)}$ is connected. 
\vspace{-0.2cm}
It is assumed that the total amount of port gain that each base station can provide is fixed. Without loss of generality, we assume $\sum_{k\in \mathcal{E}(i)}w_k^{(i)}=1$. Note that this differs from the formulation in \citep{goring2008embedded}, while they assume $\sum_{k\in\mathcal{E}} w_k=1$, and if written in the form of our model, $\sum_{i\in\mathcal{V}}\sum_{k\in\mathcal{E}(i)}w_k^{(i)}=1$. In other words, it implies that the budget of the port gains in the entire network is a constant and the power can be allocated differently to each node. Hence each node is not homogeneous any more. Therefore, the algebraic connectivity maximization of an edge-weighted Laplacian matrix can be formulated as the following SDP problem:
  \begin{equation}
	\begin{aligned}
	  & \underset{\lambda_2,\mu,\{w_k^{(i)}\}}{\text{maximize}}
	  & & \lambda_2\\
	  & \text{subject to}
	  & &\lambda_2I-\mu\mathbbm{1}\preceq \sum_{i\in\mathcal{V}}\sum_{k\in \mathcal{E}(i)}w_k^{(i)}E_k,\\
	  & & &\sum_{k\in \mathcal{E}(i)}w_k^{(i)}=1,\\
	  & & &w_k^{(i)}\geq 0,\quad \forall k\in \mathcal{E},\quad i\in \mathcal{V}.
	\end{aligned}
	\tag{P\textsubscript{C}}\label{eq:primal}
  \end{equation}

\vspace{-0.2cm}

  In \eqref{eq:primal}, the variable $\mu$ is used to shift the zero eigenvalue of $L_w$ with its eigenvector $\bm{1}$. When the optimal value is reached, $\lambda_2^*$ would be the smallest eigenvalue of $\sum_{i\in\mathcal{V}}\sum_{k\in\mathcal{E}} w_k^{(i)*}E_k+\mu^*\mathbbm{1}$. Since for any positive semi-definite matrix $G$, it holds that $\xi I\preceq G$, where $\xi$ is the smallest eigenvalue of $G$, we get the above constraints. 
  Moreover, since $\lambda_2^*$ is the second smallest eigenvalue of $\sum_{i\in\mathcal{V}}\sum_{k\in\mathcal{E}} w_k^{(i)*}E_k$, it is a continuous function with respect to $\{w_k^{(i)*}\}$; and $\{w_k^{(i)*}\}$ lives in a compact set, hence the optimal value of \eqref{eq:primal} can be attained. On the other hand, we can choose $\lambda_2$ small enough and $\mu$ large enough to make the first matrix inequality constraint strictly holds, hence strong duality holds for \eqref{eq:primal}.

\subsection{Problem Equivalence}
  In order to solve \eqref{eq:primal} distributedly, we consider the following problem
  \begin{equation}
	\begin{aligned}
	  & \underset{\{\lambda_2^{(i)}\},\{\mu^{(i)}\},\{w_k^{(i)}\},\{Z^{(i)}\}}{\text{maximize}}
	  & & \sum_{i\in\mathcal{V}}\lambda_2^{(i)}\\
	  & \text{subject to}
	  & &\lambda_2^{(i)}I-\mu^{(i)}\mathbbm{1}- \sum_{k\in \mathcal{E}(i)}w_k^{(i)}E_k\\
	  & & &+\sum_{j\in\mathcal{N}(i)}(Z^{(i)}-Z^{(j)})\preceq 0,\\
	  & & &\sum_{k\in \mathcal{E}(i)}w_k^{(i)}=1,\\
	  & & &w_k^{(i)}\geq 0,\: \forall k\in \mathcal{E}(i),\quad\forall i\in \mathcal{V},
	\end{aligned}
	\tag{P\textsubscript{D}}\label{eq:primal_distributed}
  \end{equation}
  where $Z^{(i)},\forall i\in\mathcal{V}$ are symmetric matrices. They can be written as $Z^{(i)}=\sum_{l=1}^{\frac{N(N+1)}{2}}z^{(i)}_lB_l$, where $z^{(i)}_l\in\mathbb{R}$ is the matrix entry and $B_l$ is the basis matrix for $\mathcal{S}^N$. Both are labeled by $l$. To be more precise, if $z_l^{(i)}$ is the entry located on $p$th row and $q$th column of $Z^{(i)}$, then $[B_l]_{pq}=[B_l]_{qp}=1$ and other entries of $B_l$ remain zeros. The purpose of the introduction of $Z^{(i)}$ is to derive the consensus condition \eqref{eq:consensus_constraint} in the KKT conditions.
  The next proposition describes the relationship between \eqref{eq:primal} and \eqref{eq:primal_distributed}. 
  
	\begin{prop}\label{prop:equivalence}
		If $\{\lambda_2^{(i)*}, \mu^{(i)*},\{w_k^{(i)*}\}, Z^{(i)*}\}$, $i\in\mathcal{V}$ solves \eqref{eq:primal_distributed}, then $\lambda_2^*=\sum_{i\in\mathcal{V}}\lambda_2^{(i)*}$, $\mu^*=\sum_{i\in\mathcal{V}}\mu^{(i)*}$, $\{w_k^{(i)*}\}$ solves \eqref{eq:primal}. 
		On the other hand, if there exists an optimal solution for \eqref{eq:primal}, then there also exists an optimal solution for \eqref{eq:primal_distributed}.
  \end{prop}
  \begin{pf}
	The KKT conditions of \eqref{eq:primal} for all $i\in\mathcal{V}$ are:
	\begin{subequations}\label{eq:KKT central}
	  \begin{align}
		&tr(\Phi^*)=1,\:\Phi^*\succeq 0,\:\varphi^{(i)*}_k\geq 0,\:tr(\mathbbm{1}\Phi^*)=0,\label{eq:central_der_KKT}\\
		&tr(E_k\Phi^*)-v^{(i)*}+\varphi_k^{(i)*}=0,\: w_k^{(i)*}\geq 0,\label{eq:central_E_k_Phi}\\
		&\varphi_k^{(i)*}w_k^{(i)*}=0,\: \forall k\in \mathcal{E}(i),\\
		&\lambda_2^*-\sum_{i\in\mathcal{V}}\sum_{k\in\mathcal{E}(i)}w_k^{(i)*}tr(E_k\Phi^*)=0,\label{eq:central_comple_slack}\\
		&\lambda_2^*I-\mu^*\mathbbm{1}- \sum_{i\in\mathcal{V}}\sum_{k\in \mathcal{E}(i)}w_k^{(i)*}E_k\preceq 0,\label{eq:central_feasible}\\
		&\sum_{k\in \mathcal{E}(i)}w_k^{(i)*}=1,
	  \end{align}
	\end{subequations}
	while the KKT conditions of \eqref{eq:primal_distributed} for all $i\in\mathcal{V}$ reads
	\begin{subequations}\label{eq:KKT distri}
	  \begin{align}
		&tr(\Phi^{(i)*})=1,\Phi^{(i)*}\succeq 0,\varphi^{(i)*}_k\geq 0,tr(\mathbbm{1}\Phi^{(i)*})=0,\label{eq:distri_der_KKT}\\
		&tr(E_k\Phi^{(i)*})-v^{(i)*}+\varphi_k^{(i)*}=0,\: w_k^{(i)*}\geq 0,\label{eq:E_k_Phi}\\
		&\varphi_k^{(i)*}w_k^{(i)*}=0, \forall k\in \mathcal{E}(i),\label{eq:comple_slack}\\
		&\lambda_2^{(i)*}-\sum_{k\in\mathcal{E}(i)}w_k^{(i)*}tr(E_k\Phi^{(i)*})\nonumber\\
		&\qquad\qquad\quad+\sum_{j\in\mathcal{N}(i)}tr[(Z^{(i)}-Z^{(j)})\Phi^{(i)*}]=0,\label{eq:comple_slack_Phi}\\
		&\lambda_2^{(i)*}I-\mu^{(i)*}\mathbbm{1}- \sum_{k\in \mathcal{E}(i)}w_k^{(i)*}E_k\nonumber\\
		&\qquad\qquad\qquad\qquad+\sum_{j\in\mathcal{N}(i)}(Z^{(i)*}-Z^{(j)*})\preceq 0,\label{eq:feasible}\\
		&\sum_{k\in \mathcal{E}(i)}w_k^{(i)*}=1,\quad\sum_{j\in\mathcal{N}(i)}\Phi^{(i)*}-\Phi^{(j)*}=0,\label{eq:consensus_constraint}
	  \end{align}
	\end{subequations}
	where $\Phi^*$, $\varphi_k^{(i)*}$ and $v^{(i)*}$ are the Lagrange multipliers corresponding to the matrix inequality constraint, inequality constraints and the equality constraint of \eqref{eq:primal} respectively. Similarly, $\Phi^{(i)*}$, $\varphi_k^{(i)*}$ and $v^{(i)*}$ are the Lagrange multipliers of \eqref{eq:primal_distributed}.

\vspace{-0.2cm}

	Since \eqref{eq:primal_distributed} is convex, the KKT conditions \eqref{eq:KKT distri} becomes necessary and sufficient conditions for optimality. Hence $\{\lambda_2^{(i)*}, \mu^{(i)*},\{w_k^{(i)*}\}, Z^{(i)*}\}$, $i\in\mathcal{V}$ solves \eqref{eq:primal_distributed}, if and only if there exist Lagrange multipliers  $\{\Phi^{(i)*}\}$, $\{v^{(i)*}\}$ and $\{\varphi^{(i)*}_k\}$ such that \eqref{eq:KKT distri} holds.
	Meanwhile, recall that $\{B_l\}$ are basis matrices for $\mathcal{S}^N$ and since $\Phi^{(i)*}\in\mathcal{S}^N,\forall i\in\mathcal{V}$, $\Phi^{(i)*}$ can be written as $\sum_l\phi_l^{(i)*}B_l$. Denote $\phi_l^*=[\phi_l^{(1)*},\cdots,\phi_l^{(N)*}]^T$ and \eqref{eq:consensus_constraint} can be written as $L\phi_l^*=0,\forall l$, where $L=\sum_{k\in\mathcal{E}}E_k$ is the unweighted Laplacian matrix. Since we assume that the graph is connected, then $\phi_l^*\in ker(L)= span(\bm{1}),\forall l$ and hence implies $\Phi^{(i)*}=\Phi^{(j)*}$ for all $i,j\in\mathcal{V}$.
	This means that \eqref{eq:central_der_KKT}-\eqref{eq:central_E_k_Phi} are the same as \eqref{eq:distri_der_KKT}-\eqref{eq:E_k_Phi}. 
	Further, by adding \eqref{eq:comple_slack_Phi} and \eqref{eq:feasible} for each node $i\in\mathcal{V}$, the terms $tr[(Z^{(i)}-Z^{(j)})\Phi^{(i)}]$ and $Z^{(i)*}-Z^{(j)*}$ are cancelled. Denote $\lambda_2^*=\sum_{i\in\mathcal{V}}\lambda_2^{(i)*}$, $\mu^*=\sum_{i\in\mathcal{V}}\mu^{(i)*}$, we get \eqref{eq:KKT central}. Since \eqref{eq:primal} is convex, then the KKT conditions are necessary and sufficient conditions for optimality. Hence the first part of the statement follows.

\vspace{-0.35cm}

	Now we show the second part of the statement. Suppose $\lambda_2^*,\mu^*,\{w_k^{(i)*}\}$ is an optimal solution to \eqref{eq:primal}, then there must exist Lagrange multipliers $\Phi^*$, $\{v^{(i)*}\}$ and $\{\varphi^{(i)*}_k\}$ such that \eqref{eq:KKT central} holds. Now choose $\{\hat{\lambda}_2^{(i)*},\hat{\mu}^{(i)*},\{\hat{w}_k^{(i)*}\},\hat{Z}^{(i)*}\}$ and Lagrange multipliers $\{\hat{\Phi}^{(i)*},\hat{v}^{(i)*},\{\hat{\varphi}_k^{(i)*}\}\}$ such that $\sum_{i\in\mathcal{V}}\hat{\lambda}_2^{(i)*}=\lambda_2^*,\sum_{i\in\mathcal{V}}\hat{\mu}^{(i)*}=\mu^*,\hat{w}_k^{(i)*}=w_k^{(i)*}$ and $\hat{\Phi}^{(i)*}=\Phi^*,\quad\hat{v}^{(i)*}=v^{(i)*},\quad \hat{\varphi}_k^{(i)*}=\varphi_k^{(i)*},\forall k\in\mathcal{E}(i),\forall i\in\mathcal{V}$. The KKT conditions \eqref{eq:distri_der_KKT}-\eqref{eq:comple_slack} and \eqref{eq:consensus_constraint} is trivially satisfied by the above construction. What remains to show is that there exists such $\{\hat{Z}^{(i)*}\}$ so that \eqref{eq:comple_slack_Phi} and \eqref{eq:feasible} are satisfied.
	
\vspace{-0.3cm}
We first show there exists $\{\hat{Z}^{(i)*}\}$ such that \eqref{eq:feasible} is satisfied. Denote $\hat{A}^{(i)}=\hat{\lambda}_2^{(i)*}I-\hat{\mu}^{(i)*}\mathbbm{1}-\sum_{k\in\mathcal{E}(i)}\hat{w}_k^{(i)*}E_k$. Since $\lambda_2^*,\mu^*,\{w_k^{(i)*}\}$ satisfy \eqref{eq:central_feasible} and $\lambda_2^*=\sum_{i\in\mathcal{V}}\hat{\lambda}_2^{(i)*}$, $\mu^*=\sum_{i\in\mathcal{V}}\hat{\mu}^{(i)*}$, $\hat{w}_k^{(i)*}=w_k^{(i)*}$, we know that $\sum_{i\in\mathcal{V}}\hat{A}^{(i)}\preceq 0$. By choosing $P^{(i)}=-\hat{A}^{(i)}$ for $i=1,\cdots,N-1$ and $P^{(N)}=\sum_{i=1}^{N-1}\hat{A}^{(i)}$, we have $\hat{A}^{(i)}+P^{(i)}= 0$ for $i=1,\cdots N-1$ and $\hat{A}^{(N)}+P^{(N)}=\sum_{i=1}^{N}\hat{A}^{(i)}\preceq 0$. Therefore, we have $\hat{A}^{(i)}+P^{(i)}\preceq 0$ for all $i\in\mathcal{V}$ and $\sum_{i\in\mathcal{V}}P^{(i)}=0$. What remains to show is that there exists $\{\hat{Z}^{(i)*}\}$, such that $P^{(i)}=\sum_{j\in\mathcal{N}(i)}(\hat{Z}^{(i)*}-\hat{Z}^{(j)*})$ for all $i\in\mathcal{V}$. Recall that $P^{(i)}=\sum_lp^{(i)}_lB_l$, $\hat{Z}^{(i)*}=\sum_l\hat{z}^{(i)*}_lB_l$ and denote $p_l=[p_l^{(1)},\cdots,p_l^{(N)}]^T$, $\hat{z}_l^*=[\hat{z}_l^{(1)*},\cdots,\hat{z}_l^{(N)*}]^T$. Since $\sum_{i\in\mathcal{V}}P^{(i)}=0$, it follows that $\bm{1}^Tp_l=0$ for all $l$. Therefore, $p_l\in ker(\bm{1}^T)=Im(L)$, where $L=\sum_{k\in\mathcal{E}}E_k$. This implies $p_l$ can be expressed as $L\hat{z}_l^*$ for some $\hat{z}_l^*$, namely, there exists $\{\hat{Z^{(i)}}\}$ such that $P^{(i)}=\sum_{j\in\mathcal{N}(j)}(\hat{Z}^{(i)*}-\hat{Z}^{(j)*})$ for all $i\in\mathcal{V}$ and hence \eqref{eq:feasible} is satisfied.
	
\vspace{-0.3cm}
Now we show that $\{\hat{\lambda}_2^{(i)*},\hat{\mu}^{(i)*},\hat{Z}^{(i)*}\}$ and $\{\hat{\Phi}^{(i)*}\}$ chosen above satisfy \eqref{eq:comple_slack_Phi}. Since $\hat{\Phi}^{(i)*}=\Phi^*$, $\hat{w}_k^{(i)*}=w_k^{(i)*}$ and $P^{(i)}=\sum_{j\in\mathcal{N}(i)}(\hat{Z}^{(i)*}-\hat{Z}^{(j)*})$, the left hand side of \eqref{eq:comple_slack_Phi} can be written as $\hat{\lambda}_2^{(i)*}-\sum_{k\in\mathcal{E}(i)}w_k^{(i)*}tr(E_k\Phi^{(i)*})+tr(P^{(i)}\Phi^*)$.
	Recall that $P^{(i)}=-\hat{A}^{(i)}=-\hat{\lambda}_2^{(i)*}I+\hat{\mu}^{(i)*}\mathbbm{1}+\sum_{k\in\mathcal{E}(i)}\hat{w}_k^{(i)*}E_k,i=1,\cdots N-1$. In view of \eqref{eq:central_der_KKT}, for $i=1,\cdots,N-1$, we have $\hat{\lambda}_2^{(i)*}-\sum_{k\in\mathcal{E}(i)}w_k^{(i)*}tr(E_k\Phi^*)-\hat{\lambda}_2^{(i)*}tr(\Phi^*)+\hat{\mu}^{(i)*}tr(\mathbbm{1}\Phi^*)+\sum_{k\in\mathcal{E}(i)}w_k^{(i)*}tr(E_k\Phi^*)=0$.
	
	\vspace{-0.3cm}
	For $i=N$, $P^{(N)}=\sum_{i=1}^{N-1}\hat{A}^{(i)}=\sum_{i=1}^{N-1}\hat{\lambda}_2^{(i)*}I-\sum_{i=1}^{N-1}\hat{\mu}^{(i)*}\mathbbm{1}-\sum_{i=1}^{N-1}\sum_{k\in\mathcal{E}(i)}\hat{w}_k^{(i)*}E_k$. In view of \eqref{eq:central_der_KKT}, \eqref{eq:central_comple_slack} and since $\sum_{i\in\mathcal{V}}\hat{\lambda}_2^{(i)*}=\lambda_2^*,\sum_{i\in\mathcal{V}}\hat{\mu}^{(i)*}=\mu^*$, we have $\hat{\lambda}_2^{(N)}-\sum_{k\in\mathcal{E}(N)}w_k^{(N)*}tr(E_k\Phi^*)+\sum_{i=1}^{N-1}\hat{\lambda}_2^{(i)*}tr(\Phi^*)-\sum_{i=1}^{N-1}\hat{\mu}^{(i)*}tr(\mathbbm{1}\Phi^*)-\sum_{i=1}^{N-1}\sum_{k\in\mathcal{E}(i)}w_k^{(i)*}tr(E_k\Phi^*)=0$.
	Therefore, by the variables construction above, the second part of the statement follows.\qed

  \end{pf}
  
  \subsection{Relaxing SDP into NP}  
  \eqref{eq:primal_distributed} can be solved distributedly by using a similar method as \citep{zhang2016consensus} when the graph is regular. 
  However, here we would like to consider general graphs, not only regular ones.
	In order to apply the projected saddle point dynamics to solve \eqref{eq:primal_distributed}, the problem needs to be relaxed into an NP first. This is because \eqref{eq:primal_distributed} is still an SDP problem, and its inequality matrices constraints would lead to positive semidefinite matrix Lagrangian multipliers $\Phi^{(i)}$. This makes it hard to apply the saddle point dynamics in \citep{cherukuri2016asymptotic} to this problem since by the definition of the projection operator $\Pi_K$, it is clear that $\Pi_K:\mathbb{R}^n\times\mathbb{R}^n\mapsto\mathbb{R}^n$. The projected saddle point dynamics is not defined on the cone of positive semidefinite matrices.

\vspace{-0.2cm}

  Now we introduce the convex function proposed by \cite{nesterov2007smoothing} which can be used to approximate the largest eigenvalue of a symmetric matrix. Given $X\in\mathcal{S}^N$, function $f_\varepsilon:\mathcal{S}^N\mapsto\mathbb{R}$ and reads $f_\varepsilon(X)=\varepsilon\ln tr(e^{X/\varepsilon})=\varepsilon \ln[\sum_{i=1}^Ne^{\lambda_i(X)/\varepsilon}]$ and its derivative with respect to $X$ reads
  \begin{equation}\label{eq:derivative f_epsilon}
	\nabla_Xf_\varepsilon(X)=[\sum_{i=1}^Ne^{\lambda_i(X)/\varepsilon}]^{-1}[\sum_{i=1}^Ne^{\lambda_i(X)/\varepsilon}u_iu_i^T],
  \end{equation}
  where $(\lambda_i(X),u_i)$ are eigen-pairs of $X$ with $\|u_i\|=1,\forall i$.
  It has been proved in \citep{nesterov2007smoothing} that 
   \begin{equation}
   		\lambda_{\max}(X)\leq f_\varepsilon(X)\leq \lambda_{\max}(X)+\varepsilon\ln N.
   		\label{eq:f_epsilon_approx}
   \end{equation}
  Hence when $\varepsilon$ is sufficiently small, $f_\varepsilon(X)\approx \lambda_{\max}(X)$. 

\vspace{-0.2cm}

	Consider the following NP problem:
  \begin{equation}
	\begin{aligned}
	  & \underset{\{\mu^{(i)}\},\{w_k^{(i)}\},\{Z^{(i)}\}}{\text{minimize}}
		& & \sum_{i\in\mathcal{V}}f_{\varepsilon_i}(X^{(i)})\\
	  & \text{subject to}
	  & & \sum_{k\in \mathcal{E}(i)}w_k^{(i)}=1\\
	  & & &w_k^{(i)}\geq 0,\quad \forall k\in \mathcal{E}(i),\quad\forall i\in\mathcal{V}
	\end{aligned}
	\tag{NP}\label{eq:NP}
  \end{equation}
  where $X^{(i)}=-\mu^{(i)}\mathbbm{1}-\sum_{k\in \mathcal{E}(i)}w_k^{(i)}E_k+\sum_{j\in\mathcal{N}(i)}(Z^{(i)}-Z^{(j)})$ to abbreviate the notation.
  It is clear that \eqref{eq:NP} is a convex problem, the Slater's condition also holds for \eqref{eq:NP} and hence strong duality holds. Therefore, KKT conditions becomes necessary and sufficient conditions for \eqref{eq:NP}.
The next proposition shows that how well the approximation could be.
		\begin{prop}\label{prop:approx}
		Suppose $\sum_{i\in\mathcal{V}}\lambda_2^{(i)*}$ is the optimal objective function value of \eqref{eq:primal_distributed}, then
		\begin{align}
			-\sum_{i\in\mathcal{V}}\lambda_2^{(i)*}\leq\sum_{i\in\mathcal{V}}f_{\varepsilon_i}(X^{(i)*})\leq-\sum_{i\in\mathcal{V}}\lambda_2^{(i)*}+\sum_{i\in\mathcal{V}}\varepsilon_i\ln N,
			\label{eq:approx}
		\end{align}
		where $X^{(i)*}=-\mu^{(i)*}\mathbbm{1}-\sum_{k\in \mathcal{E}(i)}w_k^{(i)*}E_k+\sum_{j\in\mathcal{N}(i)}(Z^{(i)*}-Z^{(j)*})$.
		Moreover, suppose $\{\hat{\mu}^{(i)*},\{\hat{w}_k^{(i)*}\},\hat{Z}^{(i)*}\}$ is the optimal solution to \eqref{eq:NP}, then
		\begin{equation}
			\begin{aligned}
			-\sum_{i\in\mathcal{V}}\lambda_2^{(i)*}&\leq \sum_{i\in\mathcal{V}}\lambda_{\max}(\hat{X}^{(i)*})\leq \sum_{i\in\mathcal{V}}f_{\varepsilon_i}(\hat{X}^{(i)*})\\
				&\leq -\sum_{i\in\mathcal{V}}\lambda_2^{(i)*}+\sum_{i\in\mathcal{V}}\varepsilon_i\ln N,
			\end{aligned}
			\label{eq:approx1}
		\end{equation}
		where $\hat{X}^{(i)*}=-\hat{\mu}^{(i)*}\mathbbm{1}-\sum_{k\in \mathcal{E}(i)}\hat{w}_k^{(i)*}E_k+\sum_{j\in\mathcal{N}(i)}(\hat{Z}^{(i)*}-\hat{Z}^{(j)*})$.
		\end{prop}
		\begin{pf}
			From the KKT condition \eqref{eq:feasible} of \eqref{eq:primal_distributed}, we know that $\lambda_{\max}(X^{(i)*})\leq -\lambda_2^{(i)*}$. 			
			This implies 
			\begin{align}\label{eq:eigen_ineq}
			\sum_{i\in\mathcal{V}}\lambda_{\max}(X^{(i)*})\leq-\sum_{i\in\mathcal{V}}\lambda_2^{(i)*}.
			\end{align}
			On the other hand, by Proposition \ref{prop:equivalence}, $\sum_{i\in\mathcal{V}}X^{(i)*}=-\sum_{i\in\mathcal{V}}\mu^{(i)*}\mathbbm{1}-\sum_{i\in\mathcal{V}}\sum_{k\in\mathcal{E}(i)}w_k^{(i)*}E_k=-\mu^*\mathbbm{1}-\sum_{i\in\mathcal{V}}\sum_{k\in\mathcal{E}(i)}w_k^{(i)*}E_k$, where $\lambda_2^*$, $\mu^*$ and $\{w_k^{(i)*}\}$ is the optimal solution to \eqref{eq:primal}. 
		We know from \eqref{eq:primal} that $-\lambda_2^*=\lambda_{\max}(-\mu^*\mathbbm{1}-\sum_{i\in\mathcal{V}}\sum_{k\in\mathcal{E}(i)}w_k^{(i)*})=\lambda_{\max}(\sum_{i\in\mathcal{V}}X^{(i)*})$ \citep{boyd2004convex}. Hence by Proposition \ref{prop:equivalence}, we have
		\begin{align}\label{eq:eigen_eq}
		-\sum_{i\in\mathcal{V}}\lambda_2^{(i)*}=-\lambda_2^*=\lambda_{\max}(\sum_{i\in\mathcal{V}}X^{(i)*}).
		\end{align}
		Then it follows from \eqref{eq:eigen_ineq} that $\sum_{i\in\mathcal{V}}\lambda_{\max}(X^{(i)*})\leq\lambda_{\max}(\sum_{i\in\mathcal{V}}X^{(i)*})$. On the other hand, by eigenvalue inequality, we know that $\sum_{i\in\mathcal{V}}\lambda_{\max}(X^{(i)*})\geq\lambda_{\max}(\sum_{i\in\mathcal{V}}X^{(i)*})$, and hence 
		\begin{align}\label{eq:eigen_switch}
		\sum_{i\in\mathcal{V}}\lambda_{\max}(X^{(i)*})=\lambda_{\max}(\sum_{i\in\mathcal{V}}X^{(i)*}).
		\end{align}				
		Since $\lambda_{\max}(X^{(i)*})\leq f_{\varepsilon_i}(X^{(i)*})\leq\lambda_{\max}(X^{(i)*})+\varepsilon_i\ln N$, we have $\sum_{i\in\mathcal{V}}\lambda_{\max}(X^{(i)*})\leq\sum_{i\in\mathcal{V}}f_{\varepsilon_i}(X^{(i)*})\leq\sum_{i\in\mathcal{V}}\lambda_{\max}(X^{(i)*})+\sum_{i\in\mathcal{V}}\varepsilon_i\ln N$. Hence it follows from \eqref{eq:eigen_eq} and \eqref{eq:eigen_switch} that \eqref{eq:approx} holds.
		
		\vspace{-0.3cm}
		On the other hand, by \eqref{eq:f_epsilon_approx}, it follows that $\lambda_{\max}(\hat{X}^{(i)*})\leq f_{\varepsilon_i}(\hat{X}^{(i)*})$ and hence $\sum_{i\in\mathcal{V}}\lambda_{\max}(\hat{X}^{(i)*})\leq \sum_{i\in\mathcal{V}}f_{\varepsilon_i}(\hat{X}^{(i)*})$.
		Note that $\{-\lambda_{\max}(\hat{X}^{(i)*}),\hat{\mu}^{(i)*},\{\hat{w}_k^{(i)*}\},\hat{Z}^{(i)*}\}$ is also a feasible solution to \eqref{eq:primal_distributed}, then this implies $-\sum_{i\in\mathcal{V}}\lambda_2^{(i)*}\leq \sum_{i\in\mathcal{V}}\lambda_{\max}(\hat{X}^{(i)*})$. In addition, since $\{\hat{\mu}^{(i)*},\{\hat{w}_k^{(i)*}\},\hat{Z}^{(i)*}\}$ is an optimal solution to \eqref{eq:NP}, it follows that $\sum_{i\in\mathcal{V}}f_{\varepsilon_i}(\hat{X}^{(i)*})\leq \sum_{i\in\mathcal{V}}f_{\varepsilon_i}(X^{(i)*})$. Moreover, using \eqref{eq:approx}, we have $-\sum_{i\in\mathcal{V}}\lambda_2^{(i)*}\leq \sum_{i\in\mathcal{V}}\lambda_{\max}(\hat{X}^{(i)*})\leq \sum_{i\in\mathcal{V}}f_{\varepsilon_i}(\hat{X}^{(i)*})\leq \sum_{i\in\mathcal{V}}f_{\varepsilon_i}(X^{(i)*})\leq -\sum_{i\in\mathcal{V}}\lambda_2^{(i)*}+\sum_{i\in\mathcal{V}}\varepsilon_i\ln N$, which proves the statement.	\qed
		\end{pf}
Proposition \ref{prop:approx} shows that one can have a good approximation on $\sum_{i\in\mathcal{V}}\lambda_2^{(i)}$ by choosing $\varepsilon_i$, $\forall i\in\mathcal{V}$ sufficiently small. Without losing generality, we choose $\varepsilon_i=\varepsilon$, $\forall i\in\mathcal{V}$. Note that \eqref{eq:NP} explains the reason of the introduction of $\mu$ in \eqref{eq:primal} instead of writing the constraint as $\lambda_2(I-\frac{1}{N}\mathbbm{1})\preceq \sum_{i\in\mathcal{V}}\sum_{k\in \mathcal{E}(i)}w_k^{(i)}E_k$ as in \citep{ghosh2006growing}. The $\lambda_2I$ term is needed for the relaxation. It is worth noticing that $\mu^*$ does not necessarily equals to $\lambda_2^*/N$ in \eqref{eq:primal}. In fact, any $(\lambda_2^*,\mu^*,\big\{w_k^{(i)*}\big\})$ such that $\mu^*\geq \lambda_2^*/N$ is an optimal solution for \eqref{eq:primal}. Also note that $f_\varepsilon(X)$ is not a strictly convex function though it is convex. Indeed, for all $0\leq\alpha\leq 1$, it holds that $\alpha f_\varepsilon(I)+(1-\alpha)f_\varepsilon(2I)=\alpha\varepsilon\ln(Ne^{\frac{1}{\varepsilon}})+(1-\alpha)\varepsilon\ln(Ne^{\frac{2}{\varepsilon}})=\varepsilon\ln N+2-\alpha=\varepsilon\ln(Ne^{\frac{2-\alpha}{\varepsilon}})=f_\varepsilon((2-\alpha)I)=f_\varepsilon(\alpha I+(1-\alpha)2I)$.
	Hence $f_\varepsilon(X)$ is not strictly convex.
  \subsection{Projected Dynamics and Numerical Examples}
  Now we apply the projected dynmaical system to solve \eqref{eq:NP}. To abbreviate the notation, denote $x=\big[x^{(1)T},\cdots,x^{(N)T}\big]^T$, where $x^{(i)}=\big[\mu^{(i)},\{w^{(i)}_k\},\{z^{(i)}_l\}\big]^T$ and $v=\big[v^{(1)},\cdots,v^{(N)}\big]^T$. 

\vspace{-0.2cm}
The projected dynamics for each agent $i$ is given by
  \begin{subequations}\label{eq:projected dyn}
	\begin{align}
	  \dot{\mu}^{(i)}&=\langle\nabla_{X^{(i)}}f_\varepsilon(X^{(i)}),\mathbbm{1}\rangle_M=-\frac{\partial\mathscr{L}_\mathcal{A}(x,v)}{\partial\mu^{(i)}},\\
	  \dot{w}_k^{(i)}&=\Pi_{\mathbb{R}_+}\big(w_k^{(i)},\langle\nabla_{X^{(i)}}f_\varepsilon(X^{(i)}),E_k\rangle_M-v^{(i)}\nonumber\\
	  &-(\sum_{p\in \mathcal{E}(i)}w_p^{(i)}-1)\big)\nonumber\\
	  &=\Pi_{\mathbb{R}_+}(w_k^{(i)},-\frac{\partial\mathscr{L}_\mathcal{A}(x,v)}{\partial w_k^{(i)}}),\quad \forall k\in\mathcal{E}(i),
	\end{align}
  \end{subequations}
  \begin{subequations}\label{eq:projected dyn1}
	\begin{align}
	  \dot{z}^{(i)}_l&=-\sum_{j\in\mathcal{N}(i)}\langle\nabla_{X^{(i)}}f_\varepsilon(X^{(i)})-\nabla_{X^{(j)}}f_\varepsilon(X^{(j)}),B_l\rangle_M\nonumber\\
	  &=-\frac{\partial\mathscr{L}_\mathcal{A}(x,v)}{\partial z_l^{(i)}},\\
	  \dot{v}^{(i)}&=\sum_{k\in \mathcal{E}(i)}w_k^{(i)}-1=\frac{\partial\mathscr{L}_\mathcal{A}(x,v)}{\partial v^{(i)}},
	\end{align}
  \end{subequations}
  where $\mathscr{L}_\mathcal{A}(x,v)=\sum_{i\in\mathcal{V}}\{f_\varepsilon(X^{(i)})+v^{(i)}(\sum_{k\in \mathcal{E}(i)}w_k^{(i)}-1)+\frac{1}{2}(\sum_{k\in \mathcal{E}(i)}w_k^{(i)}-1)^2\}$ and $\nabla_{X^{(i)}}f_\varepsilon(X^{(i)})$ is given by \eqref{eq:derivative f_epsilon}.
  Note that in \eqref{eq:projected dyn}, \eqref{eq:projected dyn1}, every agent only uses the information that belongs to its neighbours as well as to itself. The exchanging information for each agent is the gradient $\nabla_{X^{(i)}}f_\varepsilon(X^{(i)})$. We would like to remark that, although each time step each agent has to communicate with its neighbour a vector of size $\frac{N(N+1)}{2}$, this is the price to pay in order to solve the problem distributedly. This is because of the ``dense'' structure of the problem (since we do not make special assumptions on the graph topology) and the constraint that the communication network is the physical network itself. The reasons above have made it hard to decompose the problem into small scales.
  \begin{thm}
	The system \eqref{eq:projected dyn}, \eqref{eq:projected dyn1} is well-defined and the trajectory asymptotically converges to one of the saddle points of \eqref{eq:NP} for all initial values $\mu^{(i)}(0),z_l^{(i)}(0),v^{(i)}(0)$ $\in\mathbb{R}$, $w_k^{(i)}(0)\in\mathbb{R}_+$.
  \end{thm}
  \begin{pf}
	$f_\varepsilon(X)$ has a Lipschitz continuous gradient with respect to $x$ given that $X=\sum_ix_iA_i$, where all $A_i$ are symmetric matrices \citep{nesterov2007smoothing}. By Theorem 2.5 in \citep{nagurney2012projected}, for any initial value $\mu^{(i)}(0)\in\mathbb{R}$, $w_k^{(i)}(0)\in\mathbb{R}_+$, $v^{(i)}(0)\in\mathbb{R}$ and $z_l^{(i)}(0)\in\mathbb{R}$, there exists a unique Carath\'eodory solution which continuously depends on the initial value. Therefore the system \eqref{eq:projected dyn}, \eqref{eq:projected dyn1} is well-defined and by Proposition \ref{thm:convergence}, the system \eqref{eq:projected dyn}, \eqref{eq:projected dyn1} asymptotically converges to one of the saddle points of \eqref{eq:NP}.\qed 
  \end{pf}
It seems that when simulating the projected dynamics \eqref{eq:projected dyn},\eqref{eq:projected dyn1}, one has to do eigenvalue decomposition on $X^{(i)}$ to compute $\nabla_{X^{(i)}} f(X^{(i)})$ at each time step. However, since the factors $e^{\lambda_i(X^{(i)})/\varepsilon}$ decrease very rapidly, the gradient numerically only depends on few largest eigenvalues and correpondant eigenvectors \citep{nesterov2007smoothing}. 
  Extreme eigenvalues will converge first in numerical methods such as Arnoldi scheme, hence one does not have to do the entire eigenvalue decomposition and the numerical complexity is reduced.
  
\vspace{-0.2cm}

  By Proposition \ref{prop:approx}, one can first choose an $\varepsilon$ and get an ``optimal'' algebraic connectivity under the current choice of $\varepsilon$.
  If the approximation error $\varepsilon N\ln N$ compared to the ``optimal'' algebraic connectivity under the current choice of $\varepsilon$ is not satisfying (for example, $\varepsilon N\ln N$ is approximately 10\% of the current ``optimal'' algebraic connectivity), one can decrease $\varepsilon$ until the desired relative error is achieved.
  \begin{exmp}
	We run a simple numerical example using the graph illustrated in Fig. \ref{fig:illustration} to show that the variables do converge to the optimal solution. Using CVX, we get the optimal solution of \eqref{eq:primal}: $w_1^{(1)*}=w_2^{(3)*}=1$,  $w_1^{(2)*}=w_2^{(2)*}=0.5$ and $\lambda_2^*=1.5$. Forward Euler method is used to discretize \eqref{eq:projected dyn}, \eqref{eq:projected dyn1} and we choose $\varepsilon=0.01$, time step size $\Delta t=0.01$.
\begin{figure}[!htpb]
\centering  
\begin{minipage}[b]{0.3\textwidth} 
\includegraphics[width=1\textwidth]{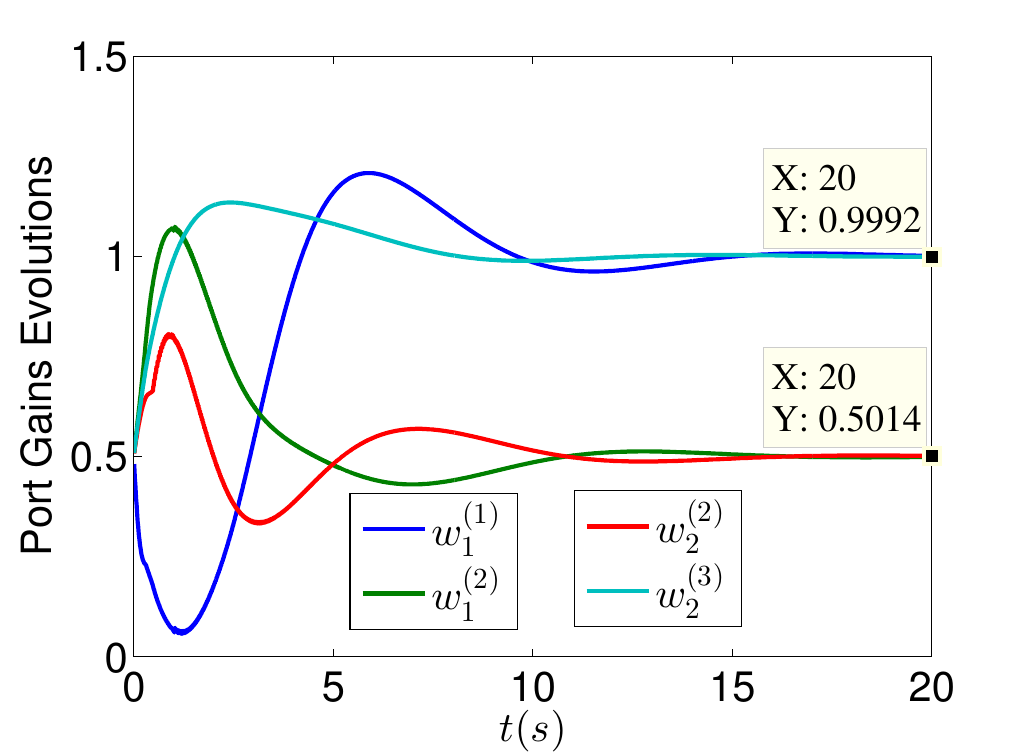} 
\end{minipage}
\caption{The evolution of the edge weights.}
\label{fig:edge weight_evolution}
\end{figure}
	Fig. \ref{fig:edge weight_evolution} shows that the edge weights converge to the optimal solution.
	\qed
  \end{exmp}
  \begin{exmp}\label{exmp:epsilon_conv_graph_topology}
	Now we consider a more complicated graph generated by ten nodes. 
	Forward Euler method is also used to discretize \eqref{eq:projected dyn}, \eqref{eq:projected dyn1} and we choose different $\varepsilon$ and time step size $\Delta t$ to illustrate the effect of $\varepsilon$ on $\Delta t$. According to \citep{nesterov2007smoothing}, the choice of $\varepsilon$ affects the Lipschitz constant of $\nabla_X f_\varepsilon(X)$ as well as the Hessian of $f_\varepsilon(X)$. The smaller $\varepsilon$ is, the bigger the Lipschitz constant of $\nabla_X f_\varepsilon(X)$ will be. Hence intuitively,  bigger Lipschitz constant of the gradient implies a smaller step size to avoid the case of moving around in the neighbourhood of optimum without converging.
	Using CVX, we know the optimal value  of \eqref{eq:primal} is 1.141. We do the simulation for $t\in[0,50]$ with different $\varepsilon$ and time step sizes. In the end, we get $\lambda_2(L_w)$ equals to 1.085, 1.091 and 1.128 when choosing $\varepsilon=10^{-2}, \Delta t=10^{-3}$, $\varepsilon=10^{-3}, \Delta t=10^{-3}$ and $\varepsilon=10^{-3}, \Delta t=10^{-4}$, respectively. As illustrated in Fig. \ref{fig:epsilon_time_step}, the algebraic connectivity in the network does not converge to the optimal value of the unrelaxed and ``centralized'' problem $(P_c)$. However, as we decrease $\varepsilon$, the limiting algebraic connectivity gets closer to the optimal value of $(P_c)$. This illustrates the relaxation effect. In addition, the evolution of $\lambda_2(L_w)$ involves a lot of oscillations when $\varepsilon=10^{-3}$ and $\Delta t=10^{-3}$, while it behaves much nicer when $\Delta t=10^{-4}$. This emperically shows that smaller $\varepsilon$ requires smaller time step length.\qed
	\begin{figure}[!htpb]
	  \centering
	  \includegraphics[width=0.35\textwidth]{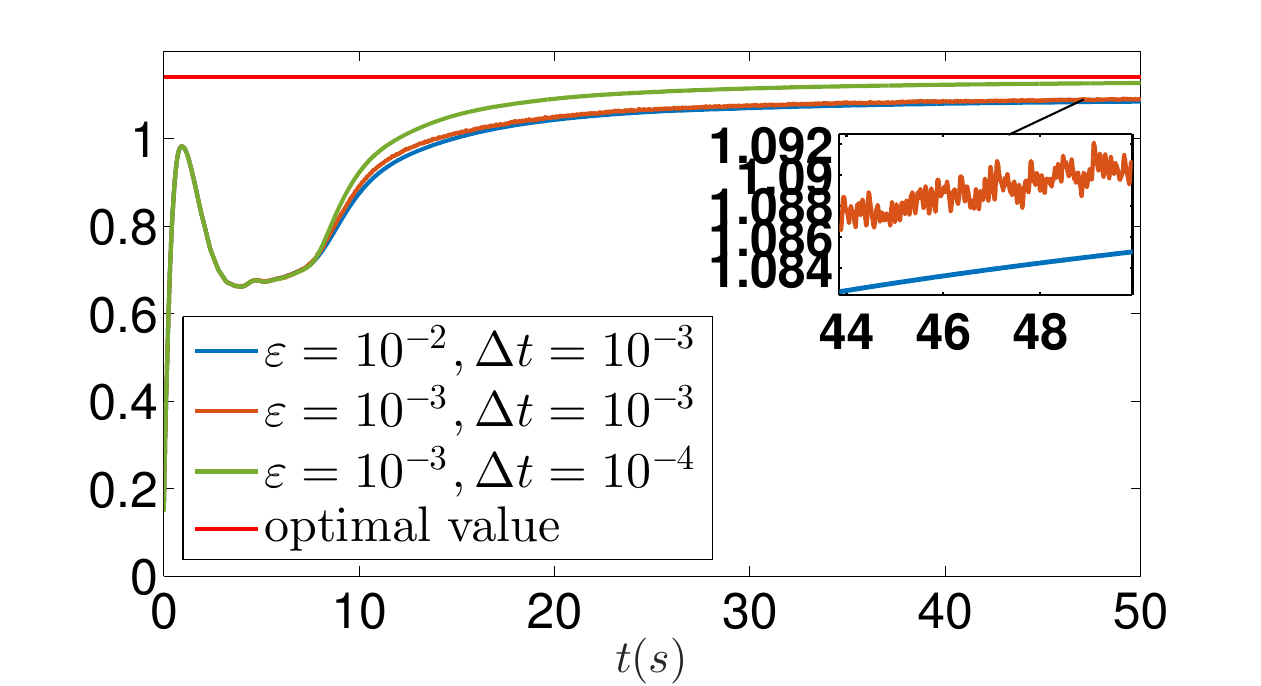}
	  \caption{The effect of relaxation and how $\varepsilon$ affects the time step size.}
	  \label{fig:epsilon_time_step}
	\end{figure}
  \end{exmp}

  Fig. \ref{fig:epsilon_time_step} illustrates that a smaller $\varepsilon$ requires a smaller step size when discretizing system \eqref{eq:projected dyn}, \eqref{eq:projected dyn1}. This means that when the number of nodes $N$ goes large, in order to get a good approximation of \eqref{eq:primal_distributed}, we need a very small $\varepsilon$ and hence it leads to a very small step size. This would result in the slow evolution of the system states per iteration and hence requires a large number of iterations to reach the equilibrium.
  
\vspace{-0.2cm}

  One practical solution to the issue above is presented as follows. We can solve the problem above by modifying $\sum_{k\in\mathcal{E}(i)}w_k^{(i)}=1$ in \eqref{eq:primal_distributed} as $\sum_{k\in\mathcal{E}(i)}w_k^{(i)}=a$, where $a>0$. We call the modified optimization problem and its relaxed nonlinear programming problem as $(P_M)$ and $(NP_M)$ respectively. By checking the optimality conditions \eqref{eq:KKT distri}, we conclude that $\{\lambda_2^{(i)*},\mu^{(i)*},\{w_k^{(i)*}\},Z^{(i)*}\}$ is the optimal solution to \eqref{eq:primal_distributed} iff $\{a\lambda_2^{(i)*},a\mu^{(i)*},a\{w_k^{(i)*}\},aZ^{(i)*}\}$ is the optimal solution to $(P_M)$ (since all the optimality conditions are linear). Using Proposition \ref{prop:approx}, we can conclude that $-\sum_{i\in\mathcal{V}}\lambda_2^{(i)*}\leq \frac{1}{a}\sum_{i\in\mathcal{V}}\lambda_{\max}(\hat{X}^{(i)*})\leq \frac{1}{a}\sum_{i\in\mathcal{V}}f_{\varepsilon_i}(\hat{X}^{(i)*})\leq -\sum_{i\in\mathcal{V}}\lambda_2^{(i)*}+\frac{1}{a}\varepsilon N\ln N$ provided that $\{\hat{\mu}^{(i)*},\{\hat{w}_k^{(i)*}\},\hat{Z}^{(i)*}\}$ is the optimal solution to $(NP_M)$. Therefore, apart from choosing $\varepsilon$ to be small, we can choose $a$ sufficiently large, solve $(NP_M)$ and divide the optimal weight realization obtained from $(NP_M)$ by $a$ to suppress the approximation error. Namely, we do not need to choose $\varepsilon$ to be too small so that the time step size does not need to be too small. Therefore, the number of iterations needed to reach the equilibrium is suppressed when $N$ goes large.
  \begin{exmp}
  Consider a graph with 30 nodes. If we do not use the methodology above, namely, $a=1$, $\varepsilon$ needs to be at least $2.1416\times 10^{-5}$ so that the relative approximation error of the optimal algebraic connectivity $\varepsilon N\ln N/\lambda_2(L_w^*)$ is within 5\%. For comparison, if we fix $\varepsilon=0.5$ first, and choose $a$ such that the relative approximation error is within 5\%. Multiple time steps have been tried and the largest ones such that the discretized systems converge are illustrated in Fig. \ref{fig:30node_a}. Same initial values and forward Euler discretization are used. The algebraic connectivity that uses the methodology mentioned above uses much fewer iterations to converge to the optimal value.
  \begin{figure}[!htpb]
	  \centering
	  \includegraphics[width=0.35\textwidth]{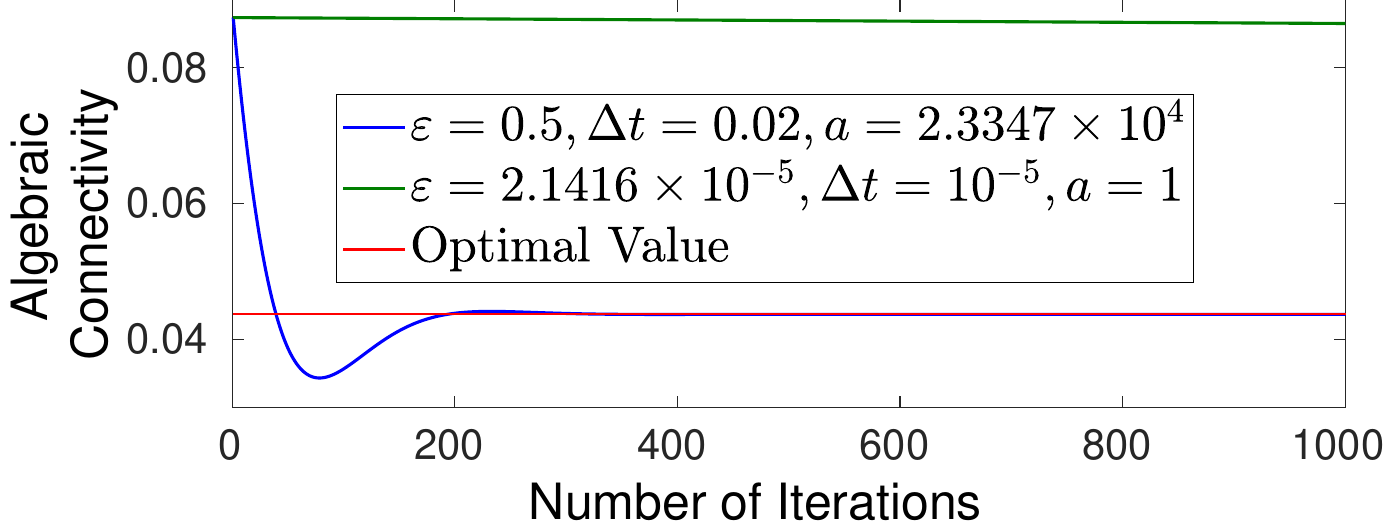}
	  \caption{The algebraic connectivities' evolution by using $a$ to suppress the number of iterations.}
	  \label{fig:30node_a}
	\end{figure}
  \end{exmp}
  \begin{exmp}
  	Consider the same graph used in Example \ref{exmp:epsilon_conv_graph_topology}. We use forward Euler for discretization. Same initial values, time step sizes and $\varepsilon$ are used ($\varepsilon=0.5$, $\Delta t=0.02$). Fig. \ref{fig:10node_a} shows how $a$ affects the approximation error. 
  	\begin{figure}[!htpb]
	  	\centering
	  	\includegraphics[width=0.45\textwidth]{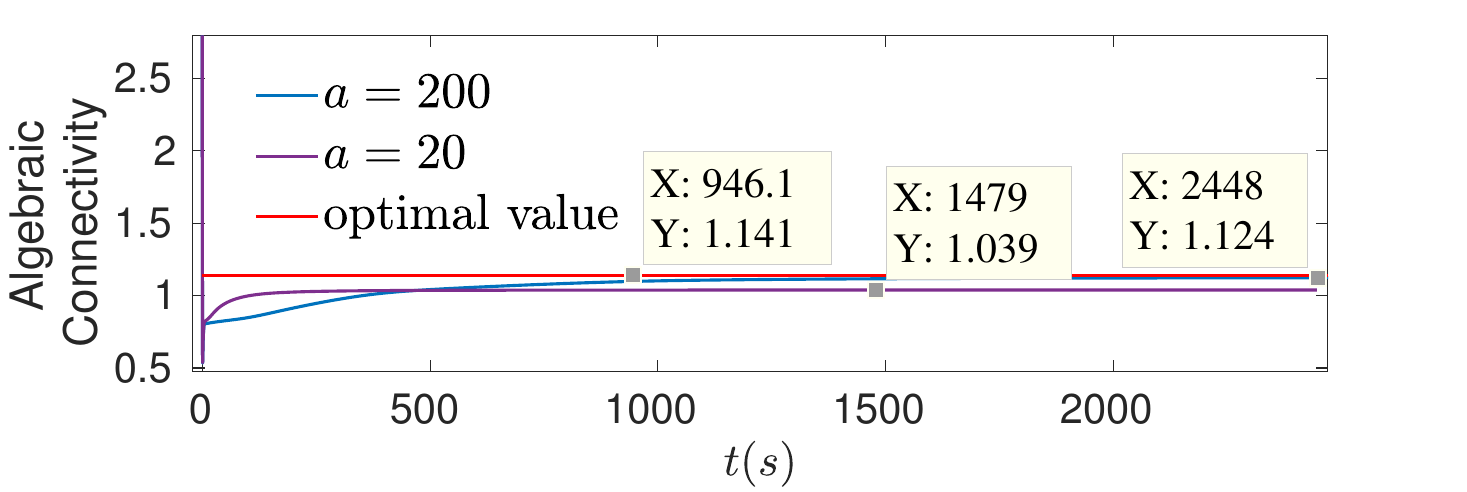}
	  	\caption{How $a$ affects the approximation error.}
	  	\label{fig:10node_a}
		\end{figure}
  \end{exmp}
  
\vspace{-0.2cm}

  The complexity per iteration for agent $i$ is $\mathcal{O}(|\mathcal{E}(i)|\cdot N^2)$. Since the right hand side of \eqref{eq:projected dyn}, \eqref{eq:projected dyn1} involves only with special matrices such as $\mathbbm{1}$, $E_k$ and $B_l$, one does not need to do matrix multiplication and hence the complexity is greatly reduced.
  \begin{exmp}
  We test the algorithm on larger scale networks and plot the ratio between running time and $|\mathcal{E}(i)|\cdot N^2$ versus $N$. Forward Euler is used for discretization. To eliminate the influence of the network topologies and number of edges on the convergence, we choose the same families of graphs and let $N$ varies. We consider the family of ring graphs and complete graphs. $\varepsilon=0.5,\Delta=0.01$ and $a$ is chosen such that the relative approximation error of the objective function is within 5\%. The iterations are terminated when the infinity norm of the right hand side of \eqref{eq:projected dyn}, \eqref{eq:projected dyn1} is smaller than $10^{-3}$. The result is shown as Fig. \ref{fig:NE}. It can be seen that the ratio between running time and $|\mathcal{E}(i)|\cdot N^2$ is approximately constant when $N$ changes.
  \begin{figure}[!htpb]
	  	\centering
	  	\includegraphics[width=0.3\textwidth]{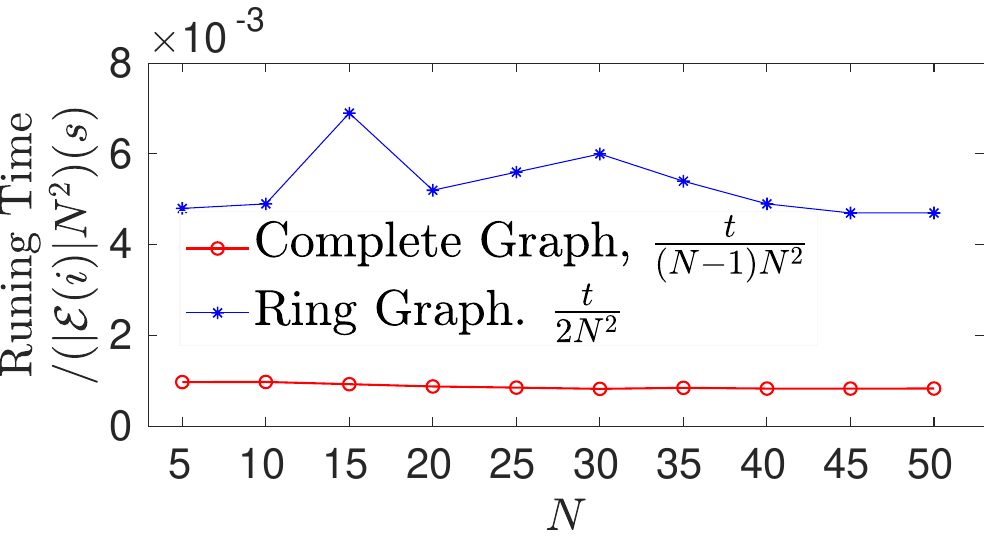}
	  	\caption{The ratio between running time and $|\mathcal{E}(i)|\cdot N^2$ versus $N$.}
	  	\label{fig:NE}
		\end{figure}
  \end{exmp}  
  
\section{Conclusion}\label{sec:conclusion}
  In this paper, a projected saddle point dynamics of augmented Lagrangian is presented to solve, not necessary strictly, convex optimization problems. As a supplement to the analysis in \citep{niederlander2016distributed}, we show that the projected saddle point dynamics converges to one of the saddle points.
  Moreover, the problem of distributedly maximizing the algebraic connectivity of an undirected communication network by optimizing the port gains of each nodes (base stations) is considered. The original SDP problem is relaxed into an NP problem and then the aforementioned projected dynamical system is applied to solve the NP. Numerical examples are used to illustrate: 1. the convergence of the edge weights to one of the optimal solutions; 2. the effect of $\varepsilon$ on the choice of time step size; 3. complexity per iteration of the algorithm. A methodology is presented so that the number of iterations needed to reach the equilibrium is suppressed. 
  
  \bibliographystyle{model5-names}
  \bibliography{projected_dyn} 

\begin{thebibliography}{27}
\expandafter\ifx\csname natexlab\endcsname\relax\def\natexlab#1{#1}\fi
\providecommand{\url}[1]{\texttt{#1}}
\providecommand{\href}[2]{#2}
\providecommand{\path}[1]{#1}
\providecommand{\DOIprefix}{doi:}
\providecommand{\ArXivprefix}{arXiv:}
\providecommand{\URLprefix}{URL: }
\providecommand{\Pubmedprefix}{pmid:}
\providecommand{\doi}[1]{\href{http://dx.doi.org/#1}{\path{#1}}}
\providecommand{\Pubmed}[1]{\href{pmid:#1}{\path{#1}}}
\providecommand{\bibinfo}[2]{#2}
\ifx\xfnm\relax \def\xfnm[#1]{\unskip,\space#1}\fi
\bibitem[{Arrow et~al.(1959)Arrow, Hurwicz, Uzawa, Chenery, Johnson, Karlin \&
  Marschak}]{arrow1959studies}
\bibinfo{author}{Arrow, K.-J.}, \bibinfo{author}{Hurwicz, L.},
  \bibinfo{author}{Uzawa, H.}, \bibinfo{author}{Chenery, H.-B.},
  \bibinfo{author}{Johnson, S.-M.}, \bibinfo{author}{Karlin, S.}, \&
  \bibinfo{author}{Marschak, T.} (\bibinfo{year}{1959}).
\newblock \bibinfo{title}{Studies in linear and non-linear programming}, .
\bibitem[{Bacciotti \& Ceragioli(2004)}]{BACCIOTTI2004841}
\bibinfo{author}{Bacciotti, A.}, \& \bibinfo{author}{Ceragioli, F.}
  (\bibinfo{year}{2004}).
\newblock \bibinfo{title}{Nonsmooth lyapunov functions and discontinuous
  carathéodory systems}.
\newblock {\it \bibinfo{journal}{IFAC Proceedings Volumes}\/},  {\it
  \bibinfo{volume}{37}\/}, \bibinfo{pages}{841 -- 845}.
\newblock \bibinfo{note}{6th IFAC Symposium on Nonlinear Control Systems 2004
  (NOLCOS 2004), Stuttgart, Germany, 2004}.
\bibitem[{Boyd \& Vandenberghe(2004)}]{boyd2004convex}
\bibinfo{author}{Boyd, S.}, \& \bibinfo{author}{Vandenberghe, L.}
  (\bibinfo{year}{2004}).
\newblock {\it \bibinfo{title}{Convex optimization}\/}.
\newblock \bibinfo{publisher}{Cambridge university press}.
\bibitem[{Brogliato et~al.(2006)Brogliato, Daniilidis, Lemar{\'e}chal \&
  Acary}]{brogliato2006equivalence}
\bibinfo{author}{Brogliato, B.}, \bibinfo{author}{Daniilidis, A.},
  \bibinfo{author}{Lemar{\'e}chal, C.}, \& \bibinfo{author}{Acary, V.}
  (\bibinfo{year}{2006}).
\newblock \bibinfo{title}{On the equivalence between complementarity systems,
  projected systems and differential inclusions}.
\newblock {\it \bibinfo{journal}{Systems \& Control Letters}\/},  {\it
  \bibinfo{volume}{55}\/}, \bibinfo{pages}{45--51}.
\bibitem[{Chatzipanagiotis et~al.(2015)Chatzipanagiotis, Dentcheva \&
  Zavlanos}]{chatzipanagiotis2015augmented}
\bibinfo{author}{Chatzipanagiotis, N.}, \bibinfo{author}{Dentcheva, D.}, \&
  \bibinfo{author}{Zavlanos, M.~M.} (\bibinfo{year}{2015}).
\newblock \bibinfo{title}{An augmented lagrangian method for distributed
  optimization}.
\newblock {\it \bibinfo{journal}{Mathematical Programming}\/},  {\it
  \bibinfo{volume}{152}\/}, \bibinfo{pages}{405--434}.
\bibitem[{Cherukuri et~al.(2015)Cherukuri, Mallada \&
  Cort{\'e}s}]{cherukuri2015convergence}
\bibinfo{author}{Cherukuri, A.}, \bibinfo{author}{Mallada, E.}, \&
  \bibinfo{author}{Cort{\'e}s, J.} (\bibinfo{year}{2015}).
\newblock \bibinfo{title}{Convergence of caratheodory solutions for primal-dual
  dynamics in constrained concave optimization}.
\newblock In {\it \bibinfo{booktitle}{SIAM conference on control and its
  applications}\/}.
\bibitem[{Cherukuri et~al.(2016)Cherukuri, Mallada \&
  Cort{\'e}s}]{cherukuri2016asymptotic}
\bibinfo{author}{Cherukuri, A.}, \bibinfo{author}{Mallada, E.}, \&
  \bibinfo{author}{Cort{\'e}s, J.} (\bibinfo{year}{2016}).
\newblock \bibinfo{title}{Asymptotic convergence of constrained primal--dual
  dynamics}.
\newblock {\it \bibinfo{journal}{Systems \& Control Letters}\/},  {\it
  \bibinfo{volume}{87}\/}, \bibinfo{pages}{10--15}.
\bibitem[{D{\"u}rr \& Ebenbauer(2011)}]{durr2011smooth}
\bibinfo{author}{D{\"u}rr, H.-B.}, \& \bibinfo{author}{Ebenbauer, C.}
  (\bibinfo{year}{2011}).
\newblock \bibinfo{title}{A smooth vector field for saddle point problems}.
\newblock In {\it \bibinfo{booktitle}{2011 50th IEEE Conference on Decision and
  Control and European Control Conference (CDC-ECC)}\/} (pp.
  \bibinfo{pages}{4654--4660}).
\newblock \bibinfo{publisher}{IEEE}.
\bibitem[{Eskelinen(2007)}]{eskelinen2007andrzej}
\bibinfo{author}{Eskelinen, P.} (\bibinfo{year}{2007}).
\newblock \bibinfo{title}{Andrzej p. ruszczy{\'n}ski: Nonlinear optimization}.
\newblock {\it \bibinfo{journal}{Mathematical Methods of Operations
  Research}\/},  {\it \bibinfo{volume}{65}\/}, \bibinfo{pages}{581--582}.
\bibitem[{Feijer \& Paganini(2010)}]{feijer2010stability}
\bibinfo{author}{Feijer, D.}, \& \bibinfo{author}{Paganini, F.}
  (\bibinfo{year}{2010}).
\newblock \bibinfo{title}{Stability of primal--dual gradient dynamics and
  applications to network optimization}.
\newblock {\it \bibinfo{journal}{Automatica}\/},  {\it \bibinfo{volume}{46}\/},
  \bibinfo{pages}{1974--1981}.
\bibitem[{Fiedler(1973)}]{fiedler1973algebraic}
\bibinfo{author}{Fiedler, M.} (\bibinfo{year}{1973}).
\newblock \bibinfo{title}{Algebraic connectivity of graphs}.
\newblock {\it \bibinfo{journal}{Czechoslovak mathematical journal}\/},  {\it
  \bibinfo{volume}{23}\/}, \bibinfo{pages}{298--305}.
\bibitem[{Ghosh \& Boyd(2006)}]{ghosh2006growing}
\bibinfo{author}{Ghosh, A.}, \& \bibinfo{author}{Boyd, S.}
  (\bibinfo{year}{2006}).
\newblock \bibinfo{title}{Growing well-connected graphs}.
\newblock In {\it \bibinfo{booktitle}{Proceedings of the 45th IEEE Conference
  on Decision and Control}\/} (pp. \bibinfo{pages}{6605--6611}).
\bibitem[{G{\"o}ring et~al.(2008)G{\"o}ring, Helmberg \&
  Wappler}]{goring2008embedded}
\bibinfo{author}{G{\"o}ring, F.}, \bibinfo{author}{Helmberg, C.}, \&
  \bibinfo{author}{Wappler, M.} (\bibinfo{year}{2008}).
\newblock \bibinfo{title}{Embedded in the shadow of the separator}.
\newblock {\it \bibinfo{journal}{SIAM Journal on Optimization}\/},  {\it
  \bibinfo{volume}{19}\/}, \bibinfo{pages}{472--501}.
\bibitem[{Khalil \& Grizzle(2002)}]{khalil1996nonlinear}
\bibinfo{author}{Khalil, H.~K.}, \& \bibinfo{author}{Grizzle, J.}
  (\bibinfo{year}{2002}).
\newblock {\it \bibinfo{title}{Nonlinear systems}\/}.
\newblock (\bibinfo{edition}{3rd} ed.).
\newblock \bibinfo{publisher}{Prentice hall New Jersey}.
\bibitem[{Kose(1956)}]{kose1956solutions}
\bibinfo{author}{Kose, T.} (\bibinfo{year}{1956}).
\newblock \bibinfo{title}{Solutions of saddle value problems by differential
  equations}.
\newblock {\it \bibinfo{journal}{Econometrica, Journal of the Econometric
  Society}\/},  (pp. \bibinfo{pages}{59--70}).
\bibitem[{Nagurney \& Zhang(2012)}]{nagurney2012projected}
\bibinfo{author}{Nagurney, A.}, \& \bibinfo{author}{Zhang, D.}
  (\bibinfo{year}{2012}).
\newblock {\it \bibinfo{title}{Projected dynamical systems and variational
  inequalities with applications}\/} volume~\bibinfo{volume}{2}.
\newblock \bibinfo{publisher}{Springer Science \& Business Media}.
\bibitem[{Nesterov(2007)}]{nesterov2007smoothing}
\bibinfo{author}{Nesterov, Y.} (\bibinfo{year}{2007}).
\newblock \bibinfo{title}{Smoothing technique and its applications in
  semidefinite optimization}.
\newblock {\it \bibinfo{journal}{Mathematical Programming}\/},  {\it
  \bibinfo{volume}{110}\/}, \bibinfo{pages}{245--259}.
\bibitem[{Niederl{\"a}nder et~al.(2016)Niederl{\"a}nder, Allg{\"o}wer \&
  Cort{\'e}s}]{niederlander2016exponentially}
\bibinfo{author}{Niederl{\"a}nder, S.~K.}, \bibinfo{author}{Allg{\"o}wer, F.},
  \& \bibinfo{author}{Cort{\'e}s, J.} (\bibinfo{year}{2016}).
\newblock \bibinfo{title}{Exponentially fast distributed coordination for
  nonsmooth convex optimization}.
\newblock In {\it \bibinfo{booktitle}{2016 IEEE 55th Conference on Decision and
  Control (CDC)}\/} (pp. \bibinfo{pages}{1036--1041}).
\newblock \bibinfo{publisher}{IEEE}.
\bibitem[{Niederl{\"a}nder \& Cort{\'e}s(2016)}]{niederlander2016distributed}
\bibinfo{author}{Niederl{\"a}nder, S.~K.}, \& \bibinfo{author}{Cort{\'e}s, J.}
  (\bibinfo{year}{2016}).
\newblock \bibinfo{title}{Distributed coordination for nonsmooth convex
  optimization via saddle-point dynamics}.
\newblock {\it \bibinfo{journal}{arXiv preprint arXiv:1606.09298}\/}, .
\bibitem[{Pakazad et~al.(2015)Pakazad, Hansson, Andersen \&
  Rantzer}]{pakazad2015distributed}
\bibinfo{author}{Pakazad, S.~K.}, \bibinfo{author}{Hansson, A.},
  \bibinfo{author}{Andersen, M.~S.}, \& \bibinfo{author}{Rantzer, A.}
  (\bibinfo{year}{2015}).
\newblock \bibinfo{title}{Distributed semidefinite programming with application
  to large-scale system analysis}.
\newblock {\it \bibinfo{journal}{arXiv preprint arXiv:1504.07755}\/}, .
\bibitem[{Schuresko \& Cort{\'e}s(2008)}]{schuresko2008distributed}
\bibinfo{author}{Schuresko, M.}, \& \bibinfo{author}{Cort{\'e}s, J.}
  (\bibinfo{year}{2008}).
\newblock \bibinfo{title}{Distributed motion constraints for algebraic
  connectivity of robotic networks}.
\newblock In {\it \bibinfo{booktitle}{2008 IEEE 47th Conference on Decision and
  Control (CDC)}\/} (pp. \bibinfo{pages}{5482--5487}).
\bibitem[{Simonetto et~al.(2013)Simonetto, Keviczky \&
  Babu{\v{s}}ka}]{simonetto2013constrained}
\bibinfo{author}{Simonetto, A.}, \bibinfo{author}{Keviczky, T.}, \&
  \bibinfo{author}{Babu{\v{s}}ka, R.} (\bibinfo{year}{2013}).
\newblock \bibinfo{title}{Constrained distributed algebraic connectivity
  maximization in robotic networks}.
\newblock {\it \bibinfo{journal}{Automatica}\/},  {\it \bibinfo{volume}{49}\/},
  \bibinfo{pages}{1348--1357}.
\bibitem[{Wang \& Elia(2011)}]{wang2011control}
\bibinfo{author}{Wang, J.}, \& \bibinfo{author}{Elia, N.}
  (\bibinfo{year}{2011}).
\newblock \bibinfo{title}{A control perspective for centralized and distributed
  convex optimization}.
\newblock In {\it \bibinfo{booktitle}{2011 50th IEEE Conference on Decision and
  Control and European Control Conference}\/} (pp.
  \bibinfo{pages}{3800--3805}).
\newblock \bibinfo{organization}{IEEE}.
\bibitem[{Yang et~al.(2010)Yang, Freeman, Gordon, Lynch, Srinivasa \&
  Sukthankar}]{yang2010decentralized}
\bibinfo{author}{Yang, P.}, \bibinfo{author}{Freeman, R.~A.},
  \bibinfo{author}{Gordon, G.~J.}, \bibinfo{author}{Lynch, K.~M.},
  \bibinfo{author}{Srinivasa, S.~S.}, \& \bibinfo{author}{Sukthankar, R.}
  (\bibinfo{year}{2010}).
\newblock \bibinfo{title}{Decentralized estimation and control of graph
  connectivity for mobile sensor networks}.
\newblock {\it \bibinfo{journal}{Automatica}\/},  {\it \bibinfo{volume}{46}\/},
  \bibinfo{pages}{390--396}.
\bibitem[{Zavlanos \& Pappas(2008)}]{zavlanos2008distributed}
\bibinfo{author}{Zavlanos, M.~M.}, \& \bibinfo{author}{Pappas, G.~J.}
  (\bibinfo{year}{2008}).
\newblock \bibinfo{title}{Distributed connectivity control of mobile networks}.
\newblock {\it \bibinfo{journal}{IEEE Transactions on Robotics}\/},  {\it
  \bibinfo{volume}{24}\/}, \bibinfo{pages}{1416--1428}.
\bibitem[{Zeng et~al.(2017)Zeng, Yi \& Hong}]{zeng2017distributed}
\bibinfo{author}{Zeng, X.}, \bibinfo{author}{Yi, P.}, \& \bibinfo{author}{Hong,
  Y.} (\bibinfo{year}{2017}).
\newblock \bibinfo{title}{Distributed continuous-time algorithm for constrained
  convex optimizations via nonsmooth analysis approach}.
\newblock {\it \bibinfo{journal}{IEEE Transactions on Automatic Control}\/},
  {\it \bibinfo{volume}{62}\/}, \bibinfo{pages}{5227--5233}.
\bibitem[{Zhang \& Hu(2016)}]{zhang2016consensus}
\bibinfo{author}{Zhang, H.}, \& \bibinfo{author}{Hu, X.}
  (\bibinfo{year}{2016}).
\newblock \bibinfo{title}{Consensus control for linear systems with optimal
  energy cost}.
\newblock {\it \bibinfo{journal}{arXiv preprint arXiv:1612.00316}\/}, .

\end{thebibliography}
\end{document}